%% file: main.tex
\documentclass[paper,11pt,nopagenum]{AAS}		

\input{preamble/packages}
\input{preamble/commands}

\input{preamble/title}

\begin{document}

\maketitle

\input{sections/0-abstract}

\section{Notation}
\begin{tabular}{ll}
$0_n$   & Matrix of zeros in $\bR{n\times n}$\\
$0_{n\times m}, 1_{n\times m}$ & Matrix of zeros and ones respectively in $\bR{n\times m}$\\
$I_n$   & Identity matrix in $\bR{n\times n}$\\
$V_{[1:n]}$ & Collection of vectors or matrices $V_k$, for $k=1,\ldots,n$\\
$v_i$ & Element $i$ of vector $v$ \\
$(u,v)$ & Concatenation of vectors $u\in\bR{n}$ and $v\in\bR{m}$ to form a vector in $\bR{n+m}$\\
$[A\,B]$ & Concatenation of matrices $A$ and $B$ with same number of rows\\
\end{tabular}

\input{sections/1-introduction}
\input{sections/2-rendezvous-problem}
\input{sections/3-scp}
\input{sections/4-pipg}
\input{sections/5-numerical-results}
\input{sections/6-conclusion}

\section{Acknowledgment}

The authors would like to thank AJ Berning for discussions on rendezvous problem formulation. We would also like to thank Purnanand Elango, Abhinav Kamath, and Yue Yu for discussions on the efficient implementation of {\pipg}. This research was supported by ONR grant N00014-20-1-2288 and Blue Origin, LLC; Government sponsorship is acknowledged.

\bibliographystyle{AAS_publication}   
\bibliography{references}   

\end{document}

%% file: preamble/packages.tex
\usepackage{bm}
\usepackage{amsmath}
\usepackage{amsfonts}
\usepackage{subfigure}
\usepackage{subcaption}
\usepackage{placeins}
\usepackage[colorlinks=false, pdfstartview=FitV, linkcolor=red, citecolor=red, urlcolor=red]{hyperref}
\usepackage{overcite}
\usepackage{algorithm}
\usepackage{algpseudocode}
\usepackage{tcolorbox}
\usepackage{tikz}
\newtcolorbox{mybox}
{
  colframe=black,
  colback=white,
  boxrule=0.75pt
}

%% file: preamble/commands.tex
\newcommand{\pipg}{\textsc{pipg}}
\newcommand{\ecos}{\textsc{ecos}}
\newcommand{\osqp}{\textsc{osqp}}

\newcommand{\D}{\mathbb{D}}
\newcommand{\R}{\mathbb{R}}

\newcommand{\range}[2]{#1\!:\!#2} 
\newcommand{\f}[2]{#1\!\left(#2\right)} 

\newcommand{\circdot}[1]{\mathring{#1}}

\newcommand{\bR}[1]{\mathbb{R}^{#1}}

%% file: preamble/title.tex
\title{Spacecraft Rendezvous Guidance via Factorization-Free Sequential Convex Programming using a First-Order Method}

\author{Govind Chari\thanks{PhD Student, Department of Aeronautics \& Astronautics, University of Washington, Seattle, WA 98195}
\ and 
{Beh\c{c}et}~{A\c{c}\i{}kme\c{s}e} \footnote{Professor, Department of Aeronautics \& Astronautics, University of Washington, Seattle, WA 98195}
}
\PaperNumber{24-162}

%% file: sections/0-abstract.tex
\begin{abstract}
We implement a fully factorization-free algorithm for nonconvex, free-final-time trajectory optimization. This algorithm is based on sequential convex programming and utilizes an inverse-free, exact discretization procedure to ensure dynamic feasibility of the converged trajectory and PIPG, a fast, first-order conic optimization algorithm as the subproblem solver. Although PIPG requires the tuning of a hyperparameter to achieve fastest convergence, we show that PIPG can be tuned to a nominal trajectory optimization problem and it is robust to variations in initial condition. We demonstrate this with a monte carlo simulation of the free-final-time rendezvous problem, using Clohessy-Wiltshire dynamics, an impulsive thrust model, and various state and control constraints including a spherical keepout zone.  
\end{abstract}

%% file: sections/1-introduction.tex
\section{Introduction}
The spacecraft rendezvous problem is concerned with finding a sequence of thrust commands and the resulting trajectory for a chaser spacecraft to bring it from some relative initial position and velocity with respect to a target spacecraft to the same final position as the target spacecraft with zero relative velocity. The trajectory and thrust commands must also satisfy a multitude of constraints such as maximum delta-v of the thrusters, maximum speed of the spacecraft, and potentially a keep-out zone.

This problem is becoming increasingly important as a number of upcoming spaceflight missions have rendezvous as a key component of their concept of operations (CONOPS). NASA is currently working on the Artemis program where they will attempt to establish long term presence on the moon \cite{chavers2020nasa}. For this project, NASA will be building the Lunar Gateway, which is a large space station in lunar orbit that will serve as a waypoint for missions from Earth to the lunar surface. SpaceX and Blue Origin have been contracted by NASA to build human lander systems (HLS), which will rendezvous with the Gateway in lunar orbit and transport astronauts to the lunar surface \cite{spacexHLS, blueHLS}. Consequently, performing a rendezvous with the Gateway is mission critical for the Artemis program.

Additionally, with companies such as Vast and Axiom Space building private space stations in low earth orbit (LEO), the number of spacecraft rendezvous performed will increase \cite{vast, axiom}. Furthermore, companies such as Starfish Space are building spacecrafts to service and deorbit satellites in LEO \cite{starfish}. All of these projects necessitates the need for reliable algorithms to solve the rendezvous problem.

A natural way of solving for trajectories subject to many constraints is formulating a nonconvex optimal control problem and solving this problem using sequential convex programming (SCP) \cite{SCPToolboxCSM2022, berning2023chance}. This is a very capable solution method and can handle difficult discrete logic constraints such as minimum impulse bit for thrusters, distance triggered plume impingement constraints, and distance triggered approach-cone constraints \cite{Malyuta2020fast, Malyuta2023fast}.

The SCP procedure requires an initial guess for the state and control trajectory. It then linearizes all nonconvex constraints and the nonlinear dynamics about this reference trajectory while keeping all convex constraints unchanged. Next, an exact discretization is performed using multiple shooting to formulate a convex subproblem which will update the reference trajectory. This subproblem can then be solved using off-the-shelf convex solvers. This process of linearization, discretization, and solving the resulting convex subproblem is repeated until convergence \cite{SCPToolboxCSM2022}. Upon convergence, the algorithm returns a dynamically feasible trajectory and control signal both of which satisfy the original constraints of the problem at the discrete nodes.

This solution procedure is agnostic to the solver used for the convex subproblem. If the subproblem is a second-order cone program (\textsc{socp}), open-source off-the-shelf solvers such as {\ecos} and {\textsc{scs}} can be used \cite{Domahidi2013ecos, scs2016, scs2021}. If the subproblem is a quadratic program (\textsc{qp}), solvers such as {\osqp} and {\textsc{piqp}} can be used in addition to the the aformentioned \textsc{socp} solvers \cite{osqp, schwan2023}. Typically the subproblem is a \textsc{socp} when constraints such as maximum speed or maximum thrust are imposed and {\ecos} is typically used as the subproblem solver \cite{Malyuta2023fast}. 

For safety critical systems such as human spaceflight, the convex solver software will need to undergo extensive testing and scrutiny in addition to verification of the underlying algorithm \cite{nasa2018software}. The {\ecos} solver is a few thousand lines of code and leverages many complex features such as sparse, permuted Cholesky factorization, Nesterov-Todd scaling, and Mehrotra's predictor-corrector \cite{Domahidi2013ecos, Andersen2011}. Fundementally understanding how this solver works requires advanced knowledge of numerical optimization which makes verification more challenging. Another downside of {\ecos} is that we cannot warmstart it. In SCP, warmstarting the solution to one subproblem based on the solution to the previous subproblem will allow quicker convex solves, but we cannot do this when using {\ecos}.

Alternatively, the Proportional-Integral Projected Gradient ({\pipg}) algorithm is four lines of pseudocode and only uses matrix-vector multiplication and closed-form projections onto simple sets such as balls and boxes \cite{yu2022extrapolated}. We can also warmstart {\pipg} which allows it to solve convex subproblems quicker. The {\pipg} algorithm has been used within the SCP framework to solve the multiphase Starship landing problem and the 6-DoF powered descent guidance problem with dual quaternions \cite{kamath2023seco, abhi2023customized}. In both cases, {\pipg} was roughly three times faster than {\ecos}. However, the downside of {\pipg} is that it requires the tuning of the hyperparameter $\omega$, the dual to primal stepsize ratio, to maximize convergence speed. 

In this paper we will solve a rendezvous problem with relevant constraints using the SCP framework with {\pipg} as a subproblem solver and show that when this hyperparameter is properly tuned to a nominal problem, it is robust to variations in problem data such as variation in initial conditions.

%% file: sections/2-rendezvous-problem.tex
\section{Rendezvous Problem Formulation}
In this section, we will formulate the free-final time, minimum control effort rendezvous problem that we will be considering in the rest of the paper.

\subsection{Dynamics}
We will begin by assuming that the target spacecraft is in a circular orbit and the chaser spacecraft is in an elliptical orbit. We can define a right handed Cartesian coordinate system whose origin is on the target spacecraft with unit vectors $\hat{i}, \hat{j}, \hat{k}$, such that $\hat{i}$ points radially outwards from the target body, $\hat{j}$ points in the direction of the velocity vector of the target body, and $\hat{k}$ points in the angular momentum direction of the target body's orbit. 

If the chaser spacecraft is sufficiently close to the target spacecraft, the equations of motion that govern the relative motion of the chaser spacecraft are the Clohessy-Wiltshire (CW) equations \cite{clohessy1960}. With the state vector defined as $x = [r^\top \; v^\top]^\top$ where $r$ is the position and $v$ is the velocity in the coordinate frame, and control input $u \in \R^3$ which is the applied force per unit mass, we can define the CW equations as follows

\begin{equation}
    \dot{x} = f(x) + Bu = \begin{bmatrix}
        v_1 \\
        v_2 \\
        v_3 \\
        3n^2r_1+2nv_2 \\
        -2nv_1 \\
        -n^2r_3
    \end{bmatrix} + \begin{bmatrix}
        0_{3 \times 3} \\
        I_3
    \end{bmatrix} u
    \label{eq:cw}
\end{equation}

where $n$ is the mean motion of the target body.

For this problem, we will model the control input as an impulse which will instantaneously change the velocity of the spacecraft. This model is valid when the amount of time taken for the spacecraft's thrusters to change the spacecraft's velocity is insignificant when compared to the total maneuver time.

If we define the change in velocity due to the thrusters for an impulsive burn as $u_b \in \R^3$ we can write the control input, $u(t)$, due to an impulse at $t_b$ as follows\cite{elango2022eigenmotion}
\begin{equation}
    u(t) = u_b\delta(t-t_b)
\end{equation}
where $\delta(t)$ is the Dirac delta function.



\subsection{Control Constraint}
We will have a control constraint that places an upper bound on the magnitude of the impulse delivered by the engine. This constraint is driven by the maximum thrust the spacecraft's propulsion system can provide. We can express this constraint as Equation \ref{eq:control-constraint} where $u_{\max}$ is the maximum delta-v the thrusters can provide.

\begin{equation}\label{eq:control-constraint}
    \|u_b\|_2 \leq u_{\max}
\end{equation}

\subsection{State Constraints}
We will consider two state constraints: a spherical keepout zone, which can easily be generalized to an elliptical keepout zone, and a maximum speed constraint. Equation \ref{eq:state-constraints} describes these constraints, where $r_c \in \R^3$ is the center of the keepout zone, $\rho_c \in \R$ is the radius of the keepout zone, and $v_{\max} \in \R$ is the maximum speed. Note that the keepout zone constraint is nonconvex.

\begin{subequations} \label{eq:state-constraints}
    \begin{align}
        &\|r(t)-r_c\|_2 \geq \rho_{c} \label{eq:nonconvex-keepout} \\
        &\|v(t)\|_2 \leq v_{\max}
    \end{align}
\end{subequations}

\subsection{Boundary Conditions}
The initial conditions for the chaser spacecraft will be some arbitrary initial position and velocity, and the terminal condition will be at the target spacecraft with zero velocity. This can be expressed as follows

\begin{subequations}
    \begin{align}
        r(0) &= r_i \\
        v(0) &= v_i \\
        r(t_f) &= 0_{3 \times 1} \\
        v(t_f) &= 0_{3 \times 1}
    \end{align}
\end{subequations}

\subsection{Continuous Time Nonconvex Optimal Control Problem}
In this problem we are minimizing control effort. We could alternatively minimize fuel consumption if we did not square the argument of the integrand, or consider a minimum time problem by minimizing $t_f$. Note that since this problem is free final time, the time-of-flight, $t_f$, is an optimization variable. We can write the full continuous time nonconvex optimal control problem as follows

\begin{mybox}\label{prob:nonconvex-ocp}
        \begin{equation*}
            \begin{aligned}
            \underset{t_{f},\,{x}(t), {u}(t)}{\text{minimize}} &
            && \int_{0}^{t_f} \|\f{u}{t}\|_2^2 \,dt \\
            \text{subject to} &
            && \forall t \in [0, t_f) \\
            \fbox{\text{CW Dynamics}} &&& {\f{\dot{{x}}}{t}} = \f{f}{{\f{{x}}{t}}, {\f{{u}}{t}}} \\
            \fbox{\text{Max delta-v}} &&& \|u_b\|_2 \le u_{\max}\\
            \fbox{\text{Keepout Zone}} &&& 
            \|{r}(t)-r_c\|_2 \geq \rho_c \hphantom{\qquad\qquad} \\
            \fbox{\text{Max Speed}}
            &&& \|{v}(t)\|_2 \leq v_{\mathrm{max}} \\
            \fbox{\text{Initial Conditions}}
            &&& {r}(0) = {r}_i \\
            &&& {v}(0) = {v}_i \\
            \fbox{\text{Terminal Conditions}}
            &&& {r}(t_f) = {0}_{3 \times 1} \\
            &&& {v}(t_f) = {0}_{3 \times 1} \\
            \end{aligned}
        \end{equation*}
\end{mybox}

%% file: sections/3-scp.tex
\section{Sequential Convex Programming}\label{sec:scp}

In this section we will outline the sequential convex programming procedure which we will use to cast the continuous time nonconvex optimal control problem in the previous section into a sequence of finite dimensional convex subproblems. We will adopt the SCP procedure outlined by Malyuta et al \cite{SCPToolboxCSM2022}, but will use time interval dilation from Kamath et al \cite{kamath2023seco} and Berning et al \cite{berning2023chance} and an inverse-free discretization procedure discussed by Kamath et al\cite{kamath2023seco}.

\subsection{Time Interval Dilation}
We apply time dilation to cast the  free final time problem as an equivalent fixed final time problem \cite{SCPToolboxCSM2022, miki2020successive}. After applying time dilation and using chain rule, our new dynamics become

\begin{equation}\label{eq:dilated-dynamics}
    \frac{dx(t(\tau))}{d\tau} = \circdot{x}(\tau) = \sigma(\tau)[f(x(\tau)) + Bu(\tau)]
\end{equation}

where $\circdot{\square}$ represents derivative with respect to $\tau \in [0, 1]$, the dilated time. We define the dilation factor as

\begin{equation}\label{eq:dilation-factor}
    \sigma(\tau) = \frac{dt}{d\tau}
\end{equation}

We can recover the time-of-flight by integrating the dilation factor as follows

\begin{equation}
    t_f = \int_0^1 \sigma(\tau)d\tau
\end{equation}

\subsection{Linearization}
After time dilation, we have a fixed final time continuous time nonconvex optimal control problem. To alleviate the nonconvexity, we linearize all nonlinear dynamics and nonconvex constraints about some reference trajectory for state, control, and dilation factor $(\bar{x}, \bar{u}, \bar{\sigma})$. To do this, we perform a first-order Taylor series expansion of \ref{eq:dilated-dynamics}.

\begin{equation}
    \circdot{x}(\tau) \approx A(\tau)x(\tau) + B(\tau)u(\tau) + S(\tau)\sigma(\tau) + c(\tau)
\end{equation}

where

\begin{subequations}
\begin{align}
A(\tau) &= \bar{\sigma} \frac{\partial f(\bar{x}(\tau), \bar{u}(\tau))}{\partial x} \\
B(\tau) &= \bar{\sigma} \frac{\partial f(\bar{x}(\tau), \bar{u}(\tau))}{\partial u} \\
S(\tau) &= f(\bar{x}, \bar{u}) \\
c(\tau) &= -A(\tau)\bar{x}(\tau) - B(\tau)\bar{u}(\tau)
\end{align}
\end{subequations}

Note that control input enters linearly into our equations of motion \ref{eq:cw}, so 

\begin{equation}\label{eq:control-jacobian}
    B(\tau) = \bar{\sigma}\begin{bmatrix}
        0_{3 \times 3} \\
        I_3
    \end{bmatrix}
\end{equation}

We must also linearize our spherical keep-out zone constraint given by equation \ref{eq:nonconvex-keepout} to obtain

\begin{equation}
    \|\bar{r}(\tau)-r_c\|_2 + \left(\frac{\bar{r}(\tau)-r_c}{\|\bar{r}(\tau)-r_c\|_2}\right)^\top (r(\tau)-\bar{r}(\tau)) \geq \rho_c
\end{equation}

\subsection{Discretization}
For computationally tractability, we must apply discretization to cast this continuous time optimal control problem as a finite dimensional optimization problem which can be solved by a computer.

We first discretize the normalized time grid $\tau \in [0, 1]$ with $K$ nodes as follows.

\begin{equation}
    0 = \tau_1 < \ldots < \tau_K = 1
\end{equation}

We will also apply a zero-order hold to the dilation factor as follows.

\begin{equation}
    \sigma(\tau) = \sigma_k \quad \forall \tau \in [\tau_k, \tau_{k+1})
\end{equation}

Since we are using an impulsive thrust model, the thrusters are only allowed to fire at the discretization nodes except the final node. Thus, $K-1$ is the maximum number of thruster firings allowed. If we used one dilation factor for the whole time grid (\emph{i.e.} $\sigma(\tau) = \sigma$), we would force all the thruster firing to be equally spaced in wall clock time. Since we apply a zero-order hold to $\sigma(\tau)$, we have a different dilation factor for each time interval, and the optimizer can choose the wall clock time between each of the burns.

We can recover the wall clock time of each interval as follows:

\begin{equation}
    t_{k+1}-t_k = \sigma_k (\tau_{k+1}-\tau_k)
\end{equation}

In order to prevent negative dilation factors we will place a lower bound on $\sigma_k$. Additionally, to prevent dilation factors that are so large that they incur extreme intersample constraint violation we will upperbound $\sigma_k$.

\begin{equation}
    \sigma_{\min} \leq \sigma_k \leq \sigma_{\max}
\end{equation}

To account for the velocity impulses, we rewrite the continuous time-dilated dynamics as a sequence of systems modeled by the unforced CW dynamics. The state and control at the $k^{\text{th}}$ node are $x_k \in \R^6$ and $u_k \in \R^3$ respectively.

\begin{subequations}
\begin{align}
        \circdot{x}(\tau) &= \sigma_kf(x(\tau)) \quad \forall \tau \in [\tau_k, \tau_{k+1}), \; k \in [1, K-1] \\
        x(\tau_k) &= x_k + (0_{3 \times 1}, u_k)
\end{align}
\end{subequations}

We can write the discretized dynamics as follows \cite{chari2024fast, kamath2023seco}

\begin{equation}
    x_{k+1} = A_k x_k + B_k u_k + S_k \sigma_k + c_k \quad \forall k \in [1, K-1]
\end{equation}

where $A_k$, $S_k$, and $c_k$ are obtained by solving the following initial value problem (IVP) on $\tau~\in~[\tau_k, \tau_{k+1}]$

\begin{subequations}
\begin{align}
\circdot{\Psi}_A &= A(\tau)\Psi_A \\
\circdot{\Psi}_S &= A(\tau)\Psi_S + S(\tau) \\
\circdot{\Psi}_c &= A(\tau)\Psi_c + c(\tau)
\end{align}
\end{subequations}

with the following initial conditions

\begin{subequations}
\begin{align}
{\Psi}_A(\tau_k) &= I_{n_x} \\
{\Psi}_S(\tau_k) &= 0_{n_x} \\
{\Psi}_c(\tau_k) &= 0_{n_x}
\end{align}
\end{subequations}

and $B_k$ is the last three columns of $A_k$, since the last three columns of $A_k$ describe how the spacecraft's velocity affects all of the spacecraft's states and our control input is an instantaneous change in velocity at the $k^\text{th}$ node.

\subsection{Virtual Terms and Trust Region}
In order to handle artificial infeasibilities due to the linearization of the nonlinear dynamics and the nonconvex spherical keep-out zone, we introduce virtual control into the dynamics and a virtual buffer term into the linearized keep-out zone. These virtual terms are slack variables which will be exactly penalized with a 1-norm in the objective function of the convex subproblem so they will be driven to zero at convergence \cite{SCPToolboxCSM2022}.

We use a first-order Taylor series expansion for all nonlinearities which is only accurate close to the reference trajectory, so we do not want the optimizer to move far from where this linearization is accurate. To encourage the optimizer to stay close to the reference trajectory, we quadratically penalize deviation from the reference trajectory in the objective function \cite{SCPToolboxCSM2022, kamath2023seco, abhi2023customized}.

After solving the convex subproblem, we will obtain a new reference trajectory and apply the same procedure of linearization, discretization, solving the convex subproblem until convergence.

\subsection{Discretized Convex Subproblem}\label{subsec:subproblem}
\begin{mybox}
        \resizebox{1.0\linewidth}{!}{$\begin{aligned}
            \underset{\sigma,\,x, u}{\text{minimize}} &
            && \sum_{k=1}^{K-1}\|u_k\|_2^2 + w_{tr}\left( \sum_{k=1}^{K}\|x_k-\bar{x}_k\|_2^2 + \sum_{k=1}^{K-1}[\|u_k-\bar{u}_k\|_2^2 + (\sigma_k-\bar{\sigma}_k)^2]\right) + w_{{vc}}\|\nu^{{c}}\|_1 + w_{{vb}}\sum_{k=1}^{K}\nu_k^b \\
            \text{subject to} &
            && \forall k \in [1, K] \\
            \fbox{\text{Discrete Dynamics}} &&& {x_{k+1}} = A_kx_k + B_ku_k + S_k\sigma_k + c_k + \nu_k^c \quad \forall k = 1,\ldots, K-1\\
            \fbox{\text{Dilation Constraints}} &&& \sigma_{\min} \leq \sigma_k \le \sigma_{\max}  \quad \forall k = 1,\ldots, K\\
            \fbox{\text{Max delta-v}} &&& \|u_k\|_2 \le u_{\max}  \quad \forall k = 1,\ldots, K-1\\
            \fbox{\text{Keepout Zone}} &&& 
            \|\bar{r}_k-r_c\|_2 + \left(\frac{\bar{r}_k-r_c}{\|\bar{r}_k-r_c\|_2}\right)^\top (r_k-\bar{r}_k) + \nu_k^b \geq \rho_c  \quad \forall k = 1,\ldots, K\\
            &&& \nu_k^b \geq 0  \quad \forall k = 1,\ldots, K\\
            \fbox{\text{Max Speed}}
            &&& \|v_k\|_2 \leq v_{\mathrm{max}}  \quad \forall k = 1,\ldots, K\\
            \fbox{\text{Initial Conditions}}
            &&& {r}(1) = {r}_i \\
            &&& {v}(1) = {v}_i \\
            \fbox{\text{Terminal Conditions}}
            &&& {r}(K) = {0}_{3 \times 1} \\
            &&& {v}(K) = {0}_{3 \times 1} \\
            \end{aligned}$}
\end{mybox}

%% file: sections/4-pipg.tex
\section{PIPG}\label{sec:pipg}
\subsection{Overview}
To solve the convex subproblems resulting from the SCP algorithm, we will use the the proportional-integral projected gradient (\pipg) algorithm. This is a first-order, primal-dual algorithm that uses closed-form projections and proportional-integral feedback of constraint violation to solve conic optimization problems \cite{yu2022extrapolated}. This algorithm has been used very effectively to solve convex 3-DoF powered-descent guidance problems, convex subproblems for 6-DoF powered-descent guidance problems, and convex subproblems in multi-phase landing problems \cite{elango2022customized, abhi2023customized, kamath2023seco}. Since it is a first-order method, the {\pipg} algorithm only involves cheap matrix-vector multiplication and avoids expensive matrix factorization.

{\pipg} solves optimization problems of the form given by Equation \ref{eq:pipgvec}, where $\D$ is the Cartesian product of convex set that have efficient projections such as balls, boxes, second-order cones, halfspaces, and the intersection of two halfspaces \cite{Bauschke2011ConvexAA}.

\begin{subequations} \label{eq:pipgvec}
\begin{align}
    \underset{z}{\text{minimize}} \quad &\frac{1}{2} z^{\top} P z + q^\top z \\
    \text{subject to} \quad & H z - h = 0 \\
        & z \in \D
\end{align}
\end{subequations}

The {\pipg} algorithm is described in Algorithm \ref{alg:pipg}, where lines 8-11 are the main loop of the algorithm and lines 1-6 involve initialization of the algorithm and computation of the primal and dual stepsizes $\alpha$ and $\beta$ respectively. In {\pipg}, $\rho \in [1.5, 1.9]$ and $\omega \in \R_{++}$ are hyperparameters. The hyperparameter $\rho$ is used in the extrapolation step, which is akin to momentum, and the hyperparameter $\omega$ is used to chose the primal and dual step sizes. For fastest convergence of {\pipg}, $\omega$ needs to be chosen in a way that the primal and dual solutions converge at the same rate.

At first glance, lines $3$ and $4$ seem expensive, however for our problem, the matrix $P$ is diagonal, so line $3$ corresponds to taking the maximum element of $P$. To compute $\max \mathrm{spec\;} H^\top H$ we can use the power iteration algorithm, which involves cheap matrix-vector products.

\begin{algorithm}[H]
\small
\caption{\pipg}\label{alg:pipg}
    \vspace{0.25em}
    \begin{flushleft}
        \textbf{Inputs:} ${P}$, ${q}$, ${H}$, $h$, ${\D}$, $\rho$, $\omega$, $k_{\max}$\\[1ex]
        \hphantom{\textbf{Inputs:}\,} $z^{0}$, ${w}^{0}$ \Comment{primal-dual guess}
    \end{flushleft}
    \vspace{0.125em}
    \begin{algorithmic}[1]
    \State $\xi^{1} \leftarrow {z}^{0}$ \Comment{initialize primal variable}
    \State $\eta^{1} \leftarrow {w}^{0}$ \Comment{initialize dual variable}\vspace{1ex}
    \State $\lambda \leftarrow \max \mathrm{spec\;} P$ \Comment{maximum eigenvalue of $P$}\vspace{1ex}
    \State $\sigma \leftarrow \max \mathrm{spec\;} H^\top H$ \Comment{maximum eigenvalue of $H^\top H$}\vspace{1ex}

    \State $\alpha \leftarrow \frac{2}{\lambda + \sqrt{\lambda^{2} + 4\omega\sigma}}$ \Comment{primal step-size}\vspace{1ex}
    \State $\beta \leftarrow \omega \alpha$ \Comment{dual step-size}\vspace{1ex}
    \For {$k \leftarrow \range{1}{k_{\max}}$}\vspace{1ex}
    \State ${z}^{k+1} \leftarrow \pi_{\D}[\xi^{k}-\alpha\,(P\,\xi^{k} + {q} + {H}^{\top} \eta^{k})]$ \Comment{primal update}
    \State ${w}^{k+1} \leftarrow \eta^{k} + \beta\,({H}(2\,{z}^{k+1}-\xi^{k}) - h)$ \Comment{dual update}
    \State $\xi^{k+1} \leftarrow (1 - \rho)\,\xi^{k} + \rho\,{z}^{k+1}$ \Comment{extrapolate primal variable}
    \State $\eta^{k+1} \leftarrow (1 - \rho)\,\eta^{k} + \rho\,{w}^{k+1}$ \Comment{extrapolate dual variable}\vspace{1ex}
    \EndFor\vspace{1ex}
    \end{algorithmic}
    \begin{flushleft}
        \textbf{Return:} $z^{k+1}$, ${w}^{k+1}$
    \end{flushleft}
    \vspace{0.25em}
\end{algorithm}

\subsection{Scaling}
Before solving the convex subproblem with PIPG we must first scale the decision variables $x, u, $ and $\sigma$ so they all have the same order of magnitude. This step is critical for the SCP algorithm to converge in practice \cite{SCPToolboxCSM2022}. Additionally, since {\pipg} is a first-order algorithm, it is not affine invariant, and improper scaling of the decision variables can result in slow convergence. Thus, we adopt the following change of variables

\begin{subequations}
\begin{align}
    x_k &= P_x \tilde{x}_k \\
    u_k &= P_u \tilde{u}_k \\
    \sigma_k &= P_\sigma \tilde{\sigma}_k
\end{align}
\end{subequations}

where $P_x$ and $P_u$ are diagonal matrices and $P_\sigma$ is a scalar. We choose these scaling parameters such that the scaled variables $\tilde{x}_k, \tilde{u}_k, $ and $\tilde{\sigma}_k$ take a maximum value of $1$. We can rewrite the discrete dynamics as follows

\begin{equation}
    {\tilde{x}_{k+1}} = \tilde{A}_k\tilde{x}_k + \tilde{B_k} \tilde{u}_k + \tilde{S}_k \tilde{\sigma}_k + \tilde{c}_k
\end{equation}

where

\begin{subequations}
\begin{align}
    \tilde{A}_k &= P_x^{-1} A_k P_x \\
    \tilde{B}_k &= P_x^{-1} B_k P_u \\
    \tilde{S}_k &= P_x^{-1} S_k P_\sigma \\
    \tilde{c}_k &= P_x^{-1} c_k
\end{align}
\end{subequations}

We must also apply this same change of variables to our reference trajectory and constraints. 

\subsection{Parsing}
Now we must \textit{parse} the discrete convex subproblem with scaled variables in the form of equation \ref{eq:pipgvec}. To keep notation light in this section, we will drop the tilde for scaled quantities. However, when solving the convex subproblem with {\pipg} or any other solver we must use the scaled quantities. 

We will define the vectorized decision variable $z$ from equation \ref{eq:pipgvec} as

\begin{equation}
    z = (x_1, \ldots, x_K, u_1, \ldots, u_{K-1}, \sigma_1, \ldots, \sigma_{K-1}, \nu_1^c, \ldots, \nu_{K-1}^c, \Gamma_1, \ldots, \Gamma_{K-1}, \nu_1^b, \ldots, \nu_{K}^b)
\end{equation}

where $x_k \in \R^{n_x}$, $u_k \in \R^{n_u}$, $\sigma_k \in \R$, $\nu^c_k \in \R^{n_x}$,  
$\Gamma_k \in \R^{n_x}$, and $\nu_k^b \in \R$. The variables $\Gamma_k \in \R^{n_x}$ are slack variables which appear when reformulating the $\|\nu^c\|_1$ term in the objective function \cite{Boyd2004}.

We can write $P$ and $q$ which define the quadratic and linear part of our cost function as follows

\begin{equation}
    P = \begin{bmatrix}
       2w_{tr}I_{n_xK} \\
    & 2(1+w_{tr})I_{n_x(K-1)} \\
    & & 2w_{tr}I_{K-1} \\
    & & & 0_{n_x(K-1)} \\
    & & & & 0_{n_x(K-1)} \\
    & & & & & 0_K \\
    \end{bmatrix}
\end{equation}

\begin{equation}
    q = (-2w_{tr}\bar{x},-2w_{tr}\bar{u}, -2w_{tr}\bar{\sigma}, 0_{n_x(K-1) \times 1}, w_{vc}1_{nx(K-1) \times 1}, w_{vb}1_{K \times 1})
\end{equation}

We will use the matrix and affine terms $H$ and $h$ to enforce just the discrete time dynamics as follows

\begin{equation}
    H = [H_x \; H_u \; H_\sigma \; H_{\nu^c} \; H_{\Gamma} \; H_{\nu^b}]
\end{equation}

\begin{equation}
    h = (-c_1, \ldots, -c_{K-1})
\end{equation}

where

\begin{subequations}
\begin{align}
    H_x &= \begin{bmatrix}
    A_1 & -I_{n_x} \\
    & A_2 & -I_{n_x} \\
    & & & \ddots \\
    & & & & A_{K-1} & -I_{n_x}
    \end{bmatrix} \\
    H_u &= \begin{bmatrix}
       B_1 \\
    & B_2 \\
    & & \ddots \\
    & & & B_{K-1} \\
    \end{bmatrix} \\
    H_\sigma &= \begin{bmatrix}
       S_1 \\
    & S_2 \\
    & & \ddots \\
    & & & S_{K-1} \\
    \end{bmatrix} \\
    H_{\nu^c} &= I_{n_x(K-1)} \\
    H_{\Gamma} &= 0_{n_x(K-1)} \\
    H_{\nu^b} &= 0_{n_x(K-1)}
\end{align}
\end{subequations}

Set $\D$ contains the remaining constraints, which are all satisfied via closed form projections.

{$\begin{aligned}
            \fbox{\text{Intersection of two halfspaces}} &&& 
            \|\hat{r}_k-r_c\|_2 + \left(\frac{\hat{r}_k-r_c}{\|\hat{r}_k-r_c\|_2}\right)^\top (r_k-\hat{r}_k) + \nu_k^b \geq \rho_c \hphantom{\qquad\qquad} \\
            &&& \nu_k^b \geq 0 \\
            \fbox{\text{Intersection of two halfspaces}} &&& 
            -\Gamma_k \leq \nu^c_k \leq \Gamma_k \\
            \fbox{\text{Box}} &&& \sigma_{\min} \leq \sigma_k \le \sigma_{\max}\\
            \fbox{\text{Ball}} &&& \|u_k\|_2 \le u_{\max}\\
            \fbox{\text{Ball}}
            &&& \|v_k\|_2 \leq v_{\mathrm{max}} \\
            \fbox{\text{Singleton}}
            &&& {r}(1) = {r}_i \\
            \fbox{\text{Singleton}}
            &&& {v}(1) = {v}_i \\
            \fbox{\text{Singleton}}
            &&& {r}(K) = {0}_{3 \times 1} \\
            \fbox{\text{Singleton}}
            &&& {v}(K) = {0}_{3 \times 1} \\
            \end{aligned}$}

Now that we have our scaled discrete time convex subproblem in the form of equation \ref{eq:pipgvec}, we can solve it with algorithm \ref{alg:pipg}. 
\subsection{Warmstarting}
In the SCP framework we solve a sequence of convex subproblems. We can use the primal-dual solution to the previous subproblem as an initial guess to warmstart {\pipg}. This enables {\pipg} to solve the subproblem in many fewer iterations than if no initial guess were provided. This technique is used to increase convergence speed of many first-order convex solvers such as {\osqp} \cite{osqp}.

%% file: sections/5-numerical-results.tex
\section{Numerical Results}

In this section we provide simulation results of the SCP algorithm described using {\pipg} to solve the recursively generated convex subproblems. We also provide a speed comparison between using {\ecos} and {\pipg} as the subproblem solver and monte carlo results to demonstrate the robustness of the SCP-{\pipg} framework to variations in initial conditions.

For all simulations, the SCP convergence criteria was the trust region radius dropping below $10^{-3}$, the virtual control $1$-norm dropping below $10^{-6}$, and the virtual buffer $1$-norm dropping below $10^{-6}$. When {\pipg} was used as a subproblem solver, it was run for $100$ iterations without a stopping criteria in place. This results in inexact solves of the subproblem, however the overall SCP algorithm converges and returns a dynamically feasible trajectory. We verify this with single shooting where we simulate the returned control signal through the original nonlinear dynamics and compare the single shooting terminal state to the desired terminal state.

\subsection{Nominal Case}
We first present the nominal rendezvous problem that we will solve using SCP with {\ecos} and {\pipg} as subproblem solvers. The rendezvous problem data and the SCP and {\pipg} hyperparameters are given in the two tables below.

    \begin{table}[h]
    \parbox{.45\linewidth}{
    \centering
    \begin{tabular}{|lll|}
    \hline
    \multicolumn{1}{|l|}{Parameter} & \multicolumn{1}{l|}{Value}              & Units                 \\ \hline
    \multicolumn{1}{|l|}{$n$}         & \multicolumn{1}{l|}{0.00113}            & $\text{s}^{-1}$ \\ \hline
    \multicolumn{1}{|l|}{$r_0$}        & \multicolumn{1}{l|}{{(}150,1000,200{)}} & m                     \\ \hline
    \multicolumn{1}{|l|}{$v_0$}        & \multicolumn{1}{l|}{{(}0,0,0{)}}        & m $\text{s}^{-1}$                   \\ \hline
    \multicolumn{1}{|l|}{$v_\text{max}$}        & \multicolumn{1}{l|}{0.5}        & m $\text{s}^{-1}$                   \\ \hline
    \multicolumn{1}{|l|}{$u_\text{max}$}        & \multicolumn{1}{l|}{0.1}        & m $\text{s}^{-1}$                   \\ \hline
    \multicolumn{1}{|l|}{$r_c$}        & \multicolumn{1}{l|}{{(}0,300,0{)}}        & m                   \\ \hline
    \multicolumn{1}{|l|}{$\rho_c$}        & \multicolumn{1}{l|}{200}        & m                   \\ \hline
    \multicolumn{1}{|l|}{$\sigma_{\text{min}}$}         & \multicolumn{1}{l|}{100}            & s \\ \hline
    \multicolumn{1}{|l|}{$\sigma_{\text{max}}$}         & \multicolumn{1}{l|}{300}            & s \\ \hline

    \end{tabular}
    \caption{Problem data}
    \label{table:nominal-data}
    }
    \hfill
    \parbox{.45\linewidth}{
    \centering
    \begin{tabular}{|ll|}
    \hline
    \multicolumn{1}{|l|}{Parameter} & \multicolumn{1}{l|}{Value}                 \\ \hline
    \multicolumn{1}{|l|}{$K$}         & \multicolumn{1}{l|}{15}            \\ \hline
    \multicolumn{1}{|l|}{$w_{{tr}}$}         & \multicolumn{1}{l|}{0.005}            \\ \hline
    \multicolumn{1}{|l|}{$w_{{vc}}$}         & \multicolumn{1}{l|}{13.0}            \\ \hline
    \multicolumn{1}{|l|}{$w_{{vb}}$}         & \multicolumn{1}{l|}{0.001}            \\ \hline
    \multicolumn{1}{|l|}{$\omega$}         & \multicolumn{1}{l|}{375.0}            \\ \hline
    \multicolumn{1}{|l|}{$\rho$}         & \multicolumn{1}{l|}{1.65}            \\ \hline
    \multicolumn{1}{|l|}{$k_{\text{max}}$}         & \multicolumn{1}{l|}{100}            \\ \hline

    \end{tabular}
    \caption{SCP-PIPG hyperparameters}
    \label{table:hyperparameters}
    }
    \end{table}

\begin{figure}[h]
    \centering
    \includegraphics[width=0.75\textwidth, trim = 15cm 2cm 12cm 2cm, clip]{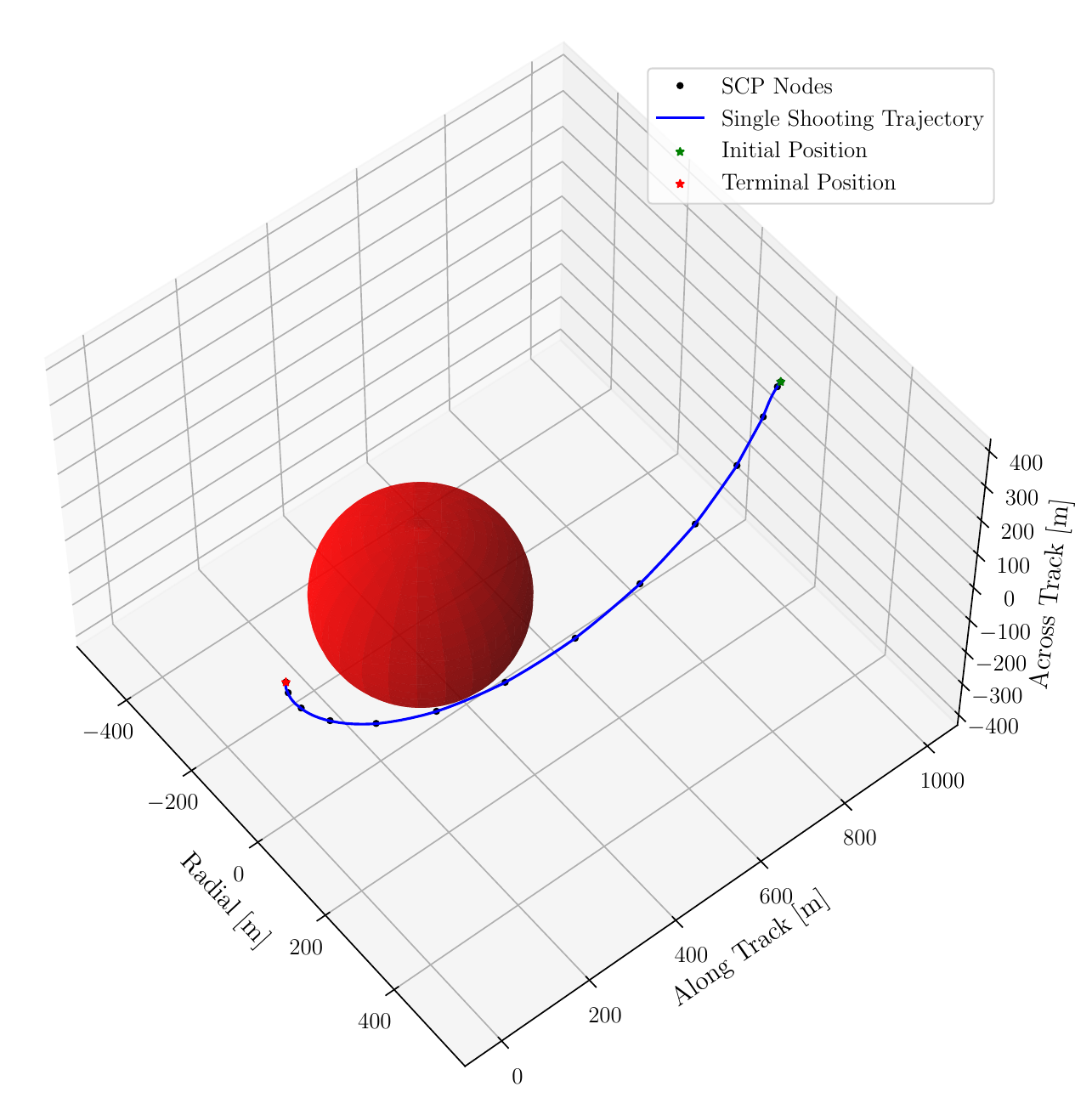}
    \caption{Rendezvous trajectory with spherical keepout}
    \label{fig:nominal-traj}
\end{figure}

\begin{figure}[h]
    \centering
    \begin{minipage}{.5\textwidth}
      \centering
    \includegraphics[width=\textwidth]{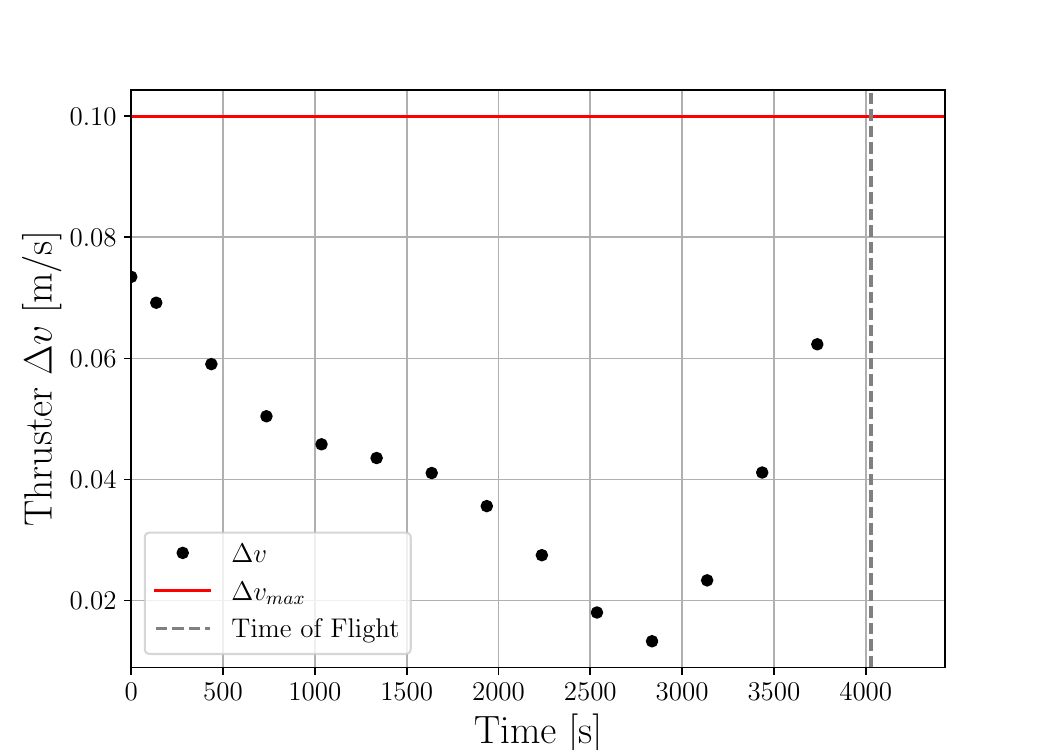}
      \caption{Thruster $\Delta v$}
      \label{fig:nominal-dv}
    \end{minipage}%
    \begin{minipage}{.5\textwidth}
      \centering
    \includegraphics[width=\textwidth]{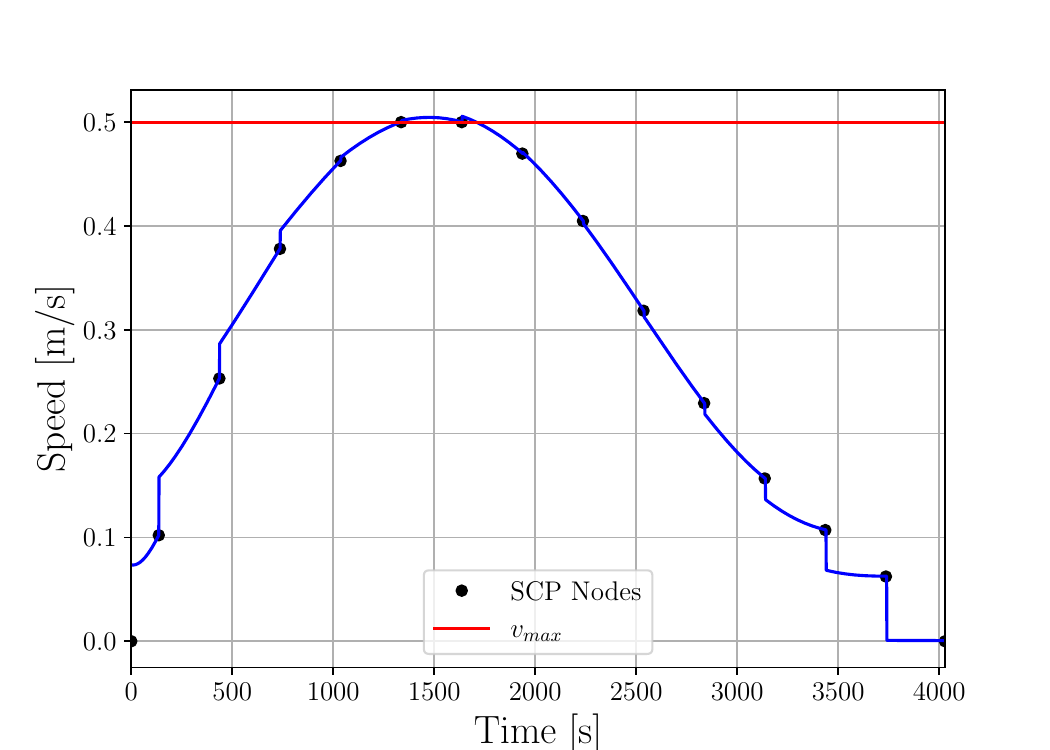}
      \caption{Vehicle speed}
      \label{fig:nominal-speed}
    \end{minipage}
\end{figure}

From Figure \ref{fig:nominal-traj} we observe that the spacecraft successfully stays out of the spherical keep-out zone, denoted by a red sphere, and the single shooting trajectory, obtained by simulating the nonlinear equations of motion with the converged control signal, passes through all the SCP nodes which indicates the dynamic feasibility of the solution \cite{miki2020successive}. We can see from Figure \ref{fig:nominal-dv} that the control signal satisfies the upper limit on $\Delta v$ we impose. We also observe that Figure \ref{fig:nominal-speed} contains discontinuities at each SCP node, which is due to the impulsive thrust model we use where the velocity of the spacecraft can be instantaneously changed by the thrusters. Furthermore, we notice that the upper speed bound is violated between the sixth and seventh SCP nodes. This intersample constraint violation is not unexpected, since we only impose constraints at the nodes \cite{dueri2017trajectory}. However, we do notice a slight violation of the speed upperbound immediately after the seventh SCP node. This is due to an interaction of how we imposed this constraint and the impulsive thrust model. We impose the state constraints at the SCP node before the impulsive thrust is applied. We can easily apply this constraint at the SCP node after the impulsive thrust is applied and we wouldn't see this speed violation.

\begin{figure}[h]
    \centering
    \includegraphics[width=0.75\textwidth]{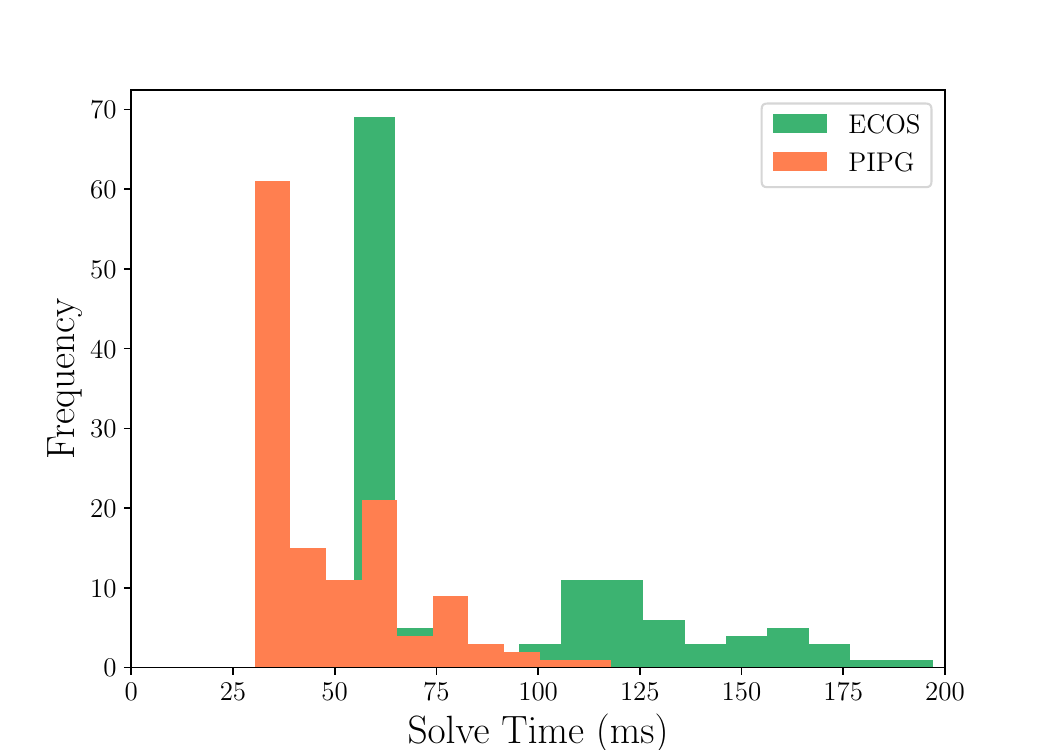}
    \caption{Runtime distribution for 128 solves of the nominal problem}
\end{figure}

In table \ref{table:nominal-stats}, we can see that when we use {\pipg} as a subproblem solver we take five more SCP iterations for the trajectory to converge due to the inexact subproblem solves. However, we obtain the same quality solution regardless of which subproblem solver we use as evidenced by the small position and velocity single shooting error. To compare the speeds of {\ecos} and {\pipg} we solved the nominal rendezvous problem 128 times with each subproblem solver and aggregated the total solve times for each SCP solve. The {\pipg} implementation is a naive implementation in Julia which performs memory allocation and solvetime and does not leverage customization, an implementation detail that exploits the structure of the subproblem and eliminates the need for sparse linear algebra, thus making solvetimes faster. Despite this and needing five more SCP iterations to converge, {\pipg} is roughly 1.82 times faster than {\ecos}.

    \begin{table}[h]
        \centering
        \begin{tabular}{|l|l|l|}
        \hline
        Statistic                            & ECOS  & PIPG \\ \hline
        SCP Iterations                       & 13    & 18   \\ \hline
        Average Subproblem Solvetime (ms)    & 90.3 & 49.5 \\ \hline
        1$\sigma$ Subproblem Solvetime (ms)    & 41.1 & 17.8 \\ \hline
        Terminal Position Single Shooting Error (m)   & 0.44     & 0.45    \\ \hline
        Terminal Velocity Single Shooting Error (m/s) & $6.4 \times 10^{-4}$     & $6.4 \times 10^{-4}$    \\ \hline
        \end{tabular}
        \caption{{\ecos} and {\pipg} solve statistics for nominal rendezvous problem}
        \label{table:nominal-stats}
    \end{table}

\FloatBarrier
\subsection{Monte Carlo Analysis}

In this section, we assess the robustness of {\pipg} as a subproblem solver within SCP to variations in initial conditions. We do this by running a monte carlo with the same problem data and SCP-{\pipg} hyperparameters given in tables \ref{table:nominal-data} and \ref{table:hyperparameters} with the exception of the initial position. We sample the initial position from the Gaussian distribution below whose mean, $r_0^{nom}$, is the nominal initial position used in the previous section, and whose standard deviation for each position component is 25 meters.

\begin{equation}
    r_0 \sim \mathcal{N}(r_0^{nom}, 25^2 I_3)    
\end{equation}

\begin{figure}[h]
    \centering
    \begin{minipage}{.5\textwidth}
      \centering
    \includegraphics[width=\textwidth]{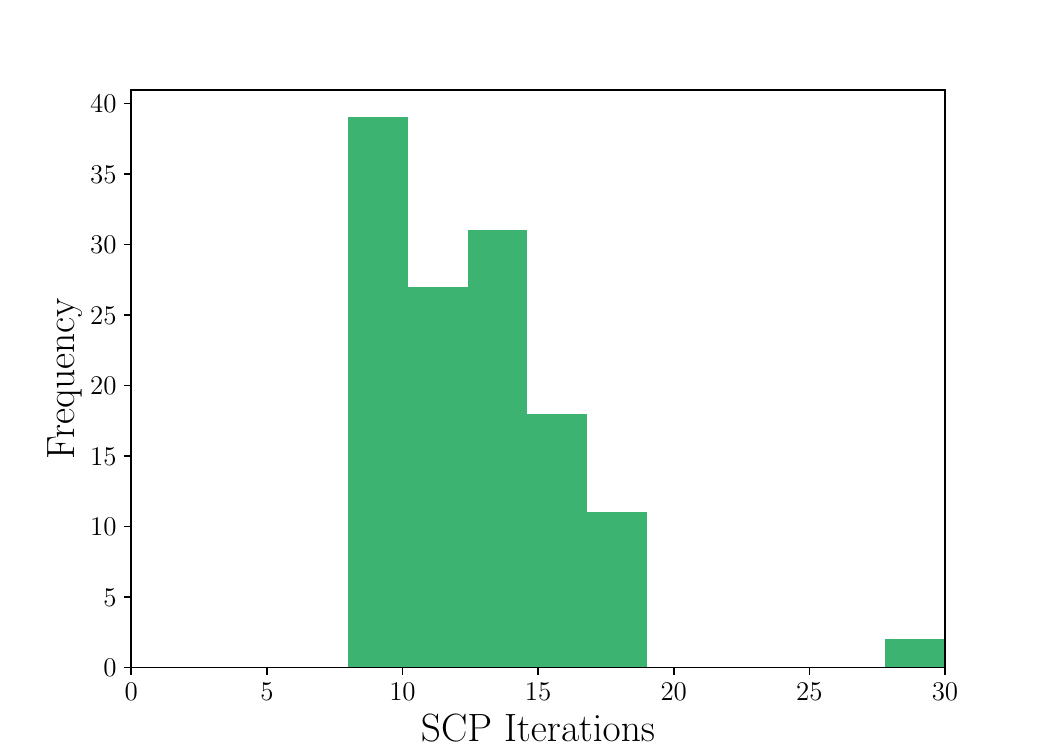}
      \caption{{\ecos} SCP iterations}
      \label{fig:ecos-scp-iters}
    \end{minipage}%
    \begin{minipage}{.5\textwidth}
      \centering
    \includegraphics[width=\textwidth]{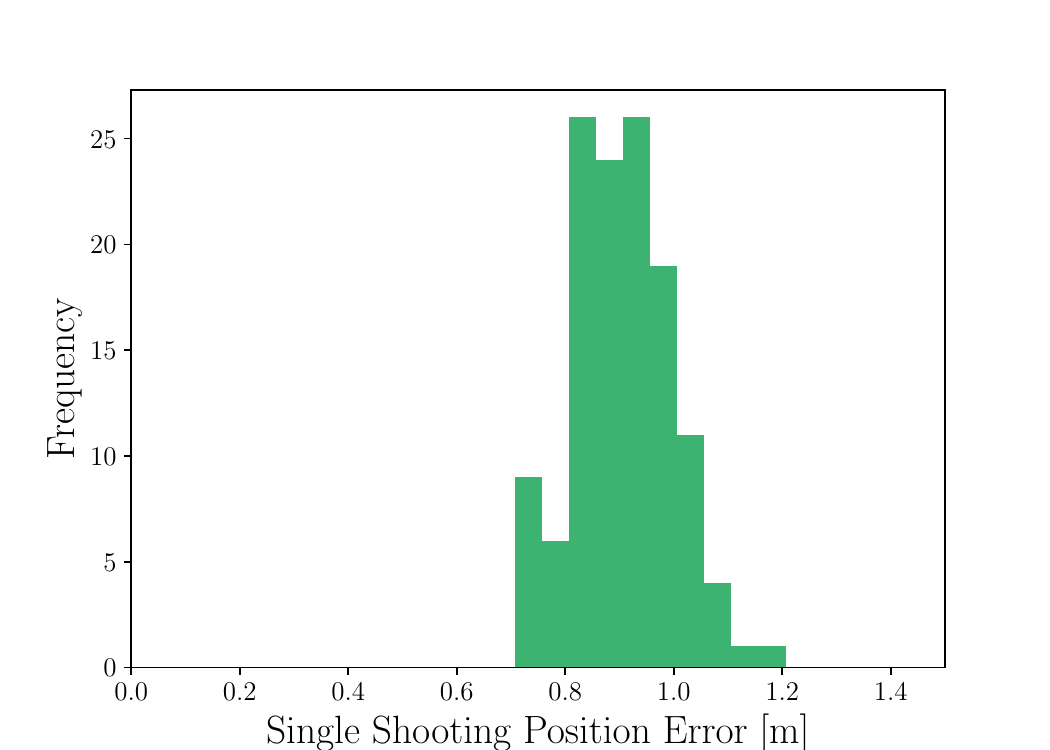}
      \caption{{\ecos} single shooting error}
      \label{fig:ecos-shooting-error}
    \end{minipage}
\end{figure}

\begin{figure}[h]
    \centering
    \begin{minipage}{.5\textwidth}
      \centering
    \includegraphics[width=\textwidth]{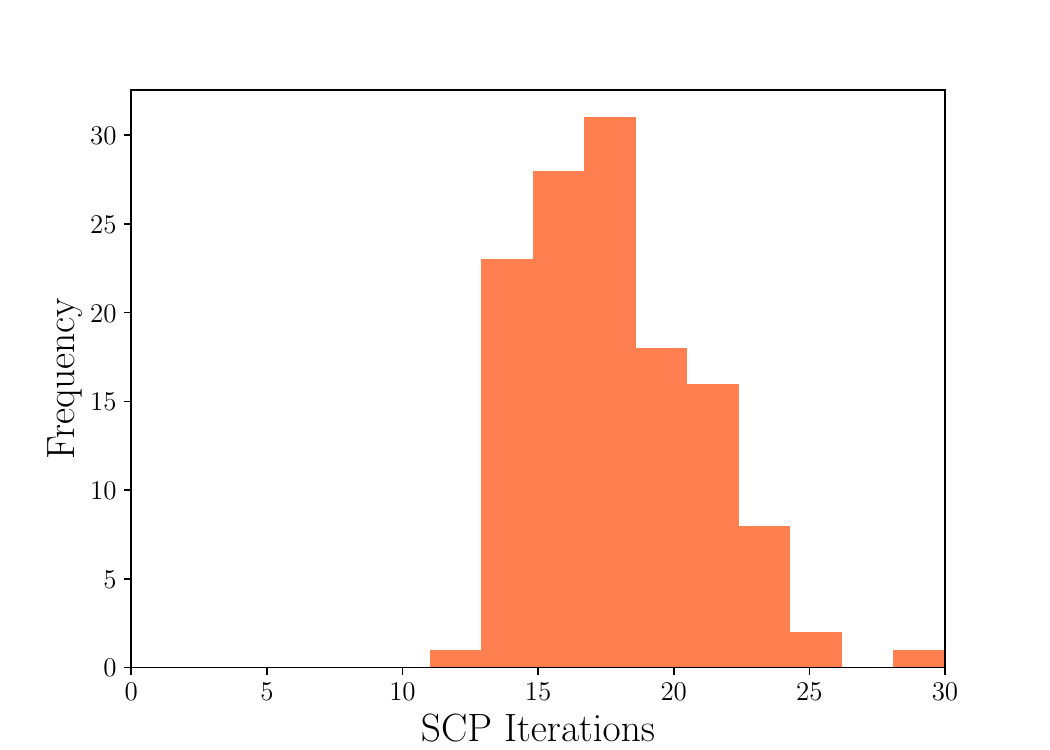}
      \caption{{\pipg} SCP iterations}
      \label{fig:pipg-scp-iters}
    \end{minipage}%
    \begin{minipage}{.5\textwidth}
      \centering
    \includegraphics[width=\textwidth]{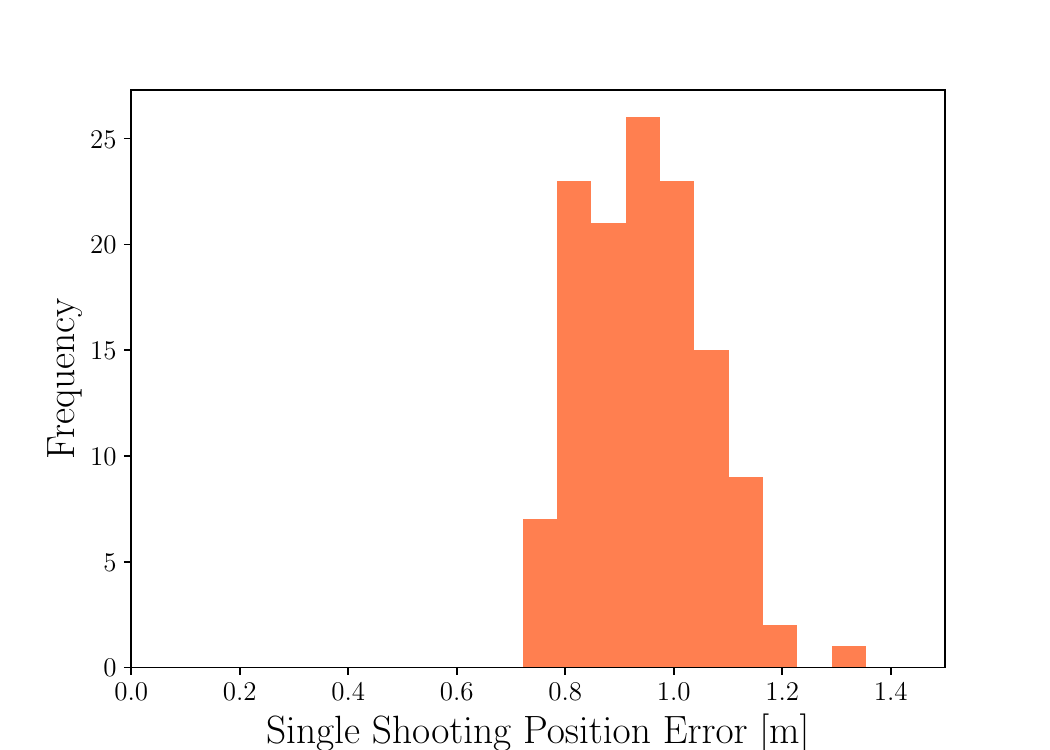}
      \caption{{\pipg} single shooting error}
      \label{fig:pipg-shooting-error}
    \end{minipage}
\end{figure}

In our monte carlo, we sampled 128 initial positions and ran SCP to generate trajectories for each of these initial conditions. From Figures \ref{fig:ecos-shooting-error} and \ref{fig:pipg-shooting-error} we can see that the single shooting position error is roughly a meter or less regardless of the subproblem solver, which indicates that our solution is dynamically feasible independent of the subproblem solver we use.

\begin{figure}[h]
    \centering
    \includegraphics[width=\textwidth, trim = 13cm 3cm 13cm 3cm, clip]{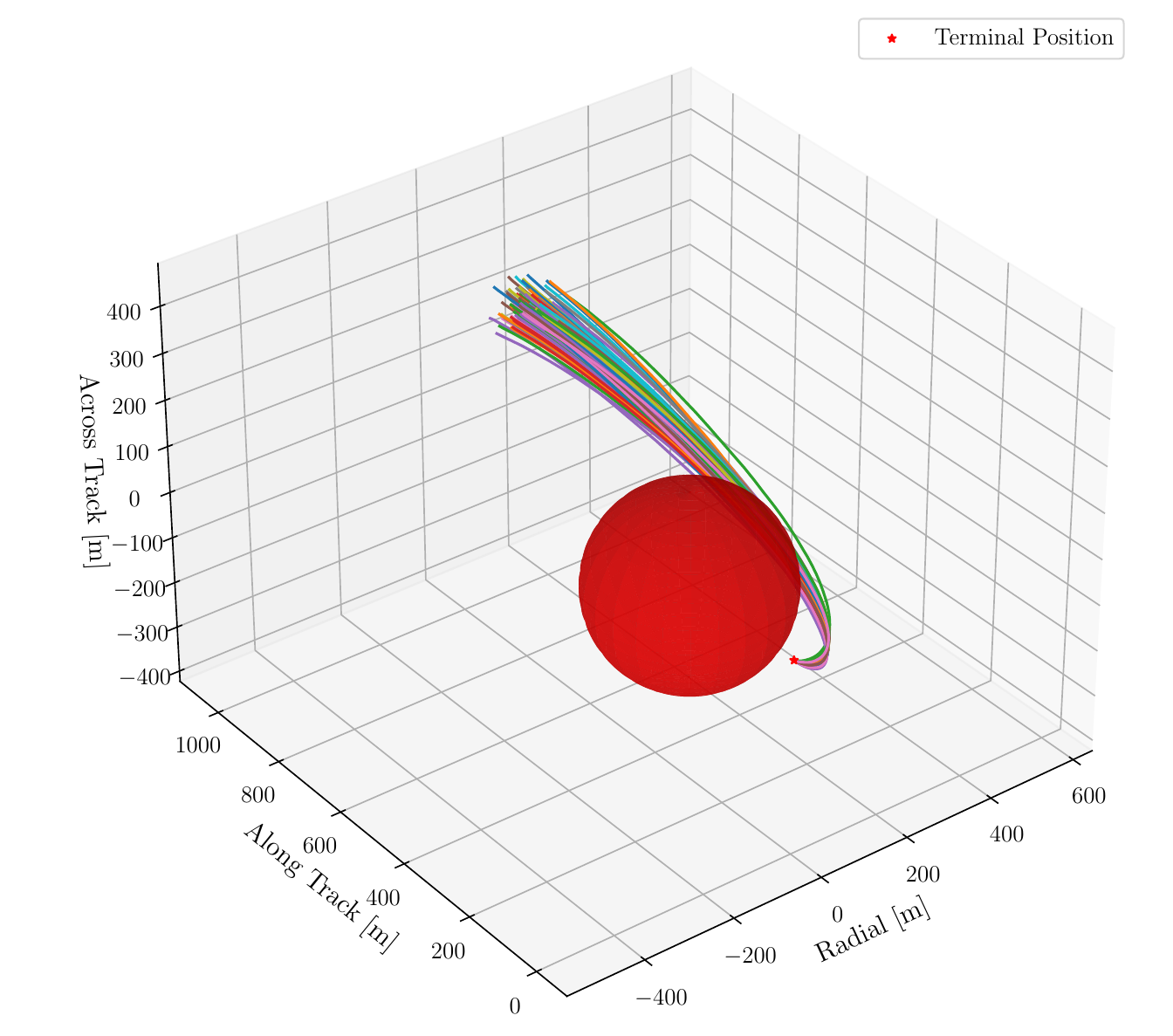}
    \caption{Monte carlo rendezvous trajectories with spherical keepout}
\end{figure}

In table \ref{table:mc-stats} we show statistics for the monte carlo simulation. We can see that for both {\ecos} and {\pipg} we have 127 converged trajectories and one trajectory which failed to converge within the 30 SCP iteration cap. From this we can see that although we tuned the {\pipg} hyperparameter, $\omega$ to the nominal rendezvous problem and ran {\pipg} for $100$ iterations for each subproblem, the solver is robust to uncertainty in the initial position. We can also see in the monte carlo that {\pipg} is roughly 1.87 times faster than {\ecos}.

    \begin{table}[h]
        \centering
        \begin{tabular}{|l|l|l|}
        \hline
        Statistic                            & ECOS  & PIPG \\ \hline
        Converged solutions               & 127/128    & 127/128   \\ \hline        
        Average SCP Iterations               & 12.6    & 17.7   \\ \hline
        1$\sigma$ SCP Iterations               & 3.5    & 3.3   \\ \hline
        Average Subproblem Solvetime (ms)    & 91.7 & 49.1 \\ \hline
        1$\sigma$ Subproblem Solvetime (ms)    & 43.0 & 17.0 \\ \hline
        Average Terminal Position Single Shooting Error (m)   & 0.91     & 0.95    \\ \hline
        1$\sigma$ Terminal Position Single Shooting Error (m)   & 0.09     & 0.11    \\ \hline
        \end{tabular}
        \caption{{\ecos} and {\pipg} solve statistics for monte carlo}
        \label{table:mc-stats}
    \end{table}

\FloatBarrier

%% file: sections/6-conclusion.tex
\section{Conclusion}
We solved the free final time nonconvex rendezvous problem using sequential convex programming in a factorization-free framework where we used the first-order conic solver {\pipg} to solve the generated sequence of convex subproblems. Although {\pipg} requires tuning of its hyperparameter, $\omega$, for fastest convergence we demonstrated that once tuned for a nominal rendezvous problem, {\pipg} is robust to variations in initial conditions. Given that it is a simpler, more explainable algorithm and we have shown it to be roughly two times faster than {\ecos}, {\pipg} is a good subproblem solver to use within the SCP framework for trajectory optimization for real-time, safety critical. For future work, we plan on developing an adaptive $\omega$ strategy based on balancing primal and dual residuals to remove the need for hyperparameter tuning in {\pipg}.